\documentclass{article}
\usepackage{amsmath, amssymb, amsfonts, amsthm, bm}

\newcommand{\Ray}[1]{\overrightarrow{\rule{0pt}{1.25ex}#1}}
\newcommand{\epsi}{\varepsilon}
\newcommand{\norm}[1]{\left\lVert#1\right\rVert}
\newcommand{\numbersystem}[1]{\mathbb{#1}}

\newcommand{\R}{\numbersystem{R}}

\newcommand{\M}{\mathcal{M}}
\newcommand{\CP}{\mathcal{P}}
\newcommand{\CS}{\mathcal{S}}
\newcommand{\abs}[1]{\lvert#1\rvert}
\newcommand{\size}[1]{\lvert#1\rvert}
\newcommand{\myangle}{\sphericalangle}

\theoremstyle{plain}
\newtheorem*{theorem}{Theorem}
\newtheorem{lemma}{Lemma}

\begin{document}

\bibliographystyle{amsplain}

\title{Low-degree minimal spanning trees in normed spaces\thanks{Research supported by a grant from an agreement between the Deutsche Forschungsgemeinschaft in Germany and the National Research Foundation in South Africa. Parts of this paper were written during a visit of the first author to the Department of Mathematical Sciences of the University of South Africa.}}
\author{Horst Martini\\
        Fakult\"at f\"ur Mathematik,\\
        Technische Universit\"at Chemnitz,\\
        D-09107 Chemnitz, Germany\\
        E-mail: \texttt{martini@mathematik.tu-chemnitz.de}
\and    Konrad J. Swanepoel\\ 
	Department of Mathematical Sciences\\
	University of South Africa \\ PO Box 392,
	Pretoria 0003 \\ South Africa \\
	E-mail: \texttt{swanekj@unisa.ac.za}}
\date{}
\maketitle

\begin{abstract}
We give a complete proof that in any finite-dimensional normed linear space a finite set of points has a minimal spanning tree in which the maximum degree is bounded above by the strict Hadwiger number of the unit ball, i.e., the largest number of unit vectors such that the distance between any two is larger than $1$.
\end{abstract}

\section{Introduction}
Let $X^d$ denote a $d$-dimensional Minkowski space, i.e., $\R^d$ equipped with a norm $\norm{\cdot}$.
Let $B=\{x:\norm{x}\leq 1\}$ be the unit ball of $X^d$, and let $B(a,r)=\{x:\norm{x-a}\leq r\}$ be the closed ball with centre $a$ and radius $r$.
For any finite set $S$ of points in $X^d$, let $\M(S)$ denote the set of minimal spanning trees (MSTs) on $S$.
As usual, a \emph{minimal spanning tree} of a set $S$ is a tree with $S$ as vertex set such that the sum of the lengths of the edges is a minimum over all trees on $S$.
Let $\Delta(T)$ denote the maximum degree of a tree $T$.
Define $\Delta^+(S)=\max\{\Delta(T):T\in\M(S)\}$ and $\Delta^-(S)=\min\{\Delta(T):T\in\M(S)\}$.
Thus $\Delta^+(S)$ is the smallest number $k$ such that all MSTs of $S$ have maximum degree at most $k$, and $\Delta^-(S)$ is the smallest number $k$ such that there exists an MST of $S$ of maximum degree at most $k$.
Now define $\Delta^+(X^d)=\max\{\Delta^+(S):S\subset X^d\}$ and $\Delta^-(X^d)=\max\{\Delta^-(S):S\subset X^d\}$.
Again, $\Delta^+(X^d)$ is the smallest number $k$ such that all MSTs of finite subsets of $X^d$ have maximum degree at most $k$, and $\Delta^-(X^d)$ is the smallest number $k$ such that any finite subset of $X^d$ has an MST with maximum degree at most $k$.

For the Euclidean plane $E^2$, for example, $\Delta^+(E^2)=6$ and $\Delta^-(E^2)=5$ \cite{PV}, while for the taxicab plane with norm $\norm{(x,y)}_1=\abs{x}+\abs{y}$ we have $\Delta^+(\R^2,\norm{\cdot}_1)=8$ \cite{HVW} and $\Delta^-(\R^2,\norm{\cdot}_1)=4$ \cite{RobSa}.

The two quantities defined above are related to two quantities from convex geometry.
Let $C$ be an arbitrary convex body in $\R^d$.
The \emph{Hadwiger number} $H(C)$ of $C$ is the largest number $k$ of translates $\{C+v_1,\dots,C+v_k\}$ of $C$ such that each $C+v_i$ touches $C$ and such that no two $C+v_i$ intersect in interior points \cite{Hadw}.
Cieslik \cite{Cieslik, Cieslik3} first proved that $\Delta^+(X^d)=H(B)$, where $B$ is the unit ball of $X^d$ (see also \cite{RobSa}).
Note that for a unit ball $B$ the value $H(B)$ equals the maximum number of unit vectors in $X^d$ such that the distance between any two is at most $1$.

The value $\Delta^-(X^d)$ also equals a similar quantity, called the MST number in \cite{RobSa}, and called the weak Hadwiger number in \cite{Swa}.
Here we propose the following more descriptive name.
The \emph{strict Hadwiger number} $H_s(C)$ of $C$ is the largest number $k$ of translates $\{C+v_1,\dots,C+v_k\}$ of $C$ such that each $C+v_i$ touches $C$ and such that any two $C+v_i$ are disjoint.
Again note that for a unit ball $B$ the value $H_s(B)$ equals the maximum number of unit vectors in $X^d$ such that the distance between any two is greater than $1$.
For the $\ell_p$ norm $\norm{(x_1,\dots,x_d)}_p=(\sum_{i=1}^d\abs{x_i}^p)^{1/p}$ with unit ball $B_p$ it was shown by Robins and Salowe \cite{RobSa} that $\Delta^-(\R^d,\norm{\cdot}_p)=H_s(B_p)$.
However, their proof uses a certain perturbation of points.
It is not immediately clear how such a perturbation is to be done.
It is the purpose of this note to clear up this point (using the Baire category theorem from point-set topology) and, simultaneously, to extend the result to \emph{arbitrary Minkowski spaces} (using planar Minkowski geometry).

\begin{theorem}
For any Minkowski space $X^d$ with unit ball $B$, \[\Delta^-(X^d)=H_s(B).\]
\end{theorem}

We note that there are many results about the Hadwiger number; see for example \cite[\S2.9, \S9.6, \S9.7]{Bor}.
We mention the following facts.
For planar convex bodies, the Hadwiger number is $8$ for parallelograms and $6$ for all other bodies \cite{MR25:1492}, for the three-dimensional octahedron $O$  (the unit ball of the $L_1$ norm in $\R^3$) it is $H(O)=18$ \cite{RobSa}, and for the $d$-cube $C_d$ (the unit ball of the $L_\infty$ norm in $\R^d$) it is easily seen that $H(C_d)=3^d-1$.
Much less is known about the strict Hadwiger number.
For planar convex bodies, the strict Hadwiger number is $4$ for parallelograms \cite{RobSa} and $5$ for all other bodies \cite{MR93e:52038}.
For the three-dimensional Euclidean ball it is easily seen that the value equals the Hadwiger number, namely $12$.
For the three-dimensional octahedron it is known that $13\leq H_s(O)\leq 14$ \cite{RobSa}, and for the $d$-cube $H_s(C_d)=2^d$ \cite{RobSa}.

\section{Proof}
There exists a set $S$ of $H_s(B)$ points on the boundary of the unit ball $B$ such that the distance between any two is greater than $1$.
Then the set $S\cup\{o\}$ has only one MST, with the origin $o$ of degree $H_s(B)$.
It follows that $\Delta^-(X^d)\geq H_s(B)$.

We now show that for any finite set $S$ in $X^d$ there exists an MST $T$ with $\Delta(T)\leq H_s(B)$.
In order to do this, we consider angles.
Any three distinct points $a,b,c\in S$ define an angle $\myangle abc$ at $b$, bounded by the rays $\Ray{ba}$ and $\Ray{bc}$.
(If $b$ is between $a$ and $c$ on the same line, we may take either half plane to be the angle.)
We define the \emph{size} of $\myangle abc$ by
\[ \size{\myangle abc} := \norm{\frac{1}{\norm{a-b}}(a-b)-\frac{1}{\norm{c-b}}(c-b)}.\]
This is the distance between the two points where the rays of the angle intersect the unit ball with centre $b$.
In Euclidean space we have that $\size{\myangle abc}=1$ if and only if the ordinary angular measure of $\myangle abc$ is $60^\circ$.
It is well-known that angles between incident edges in MSTs in Euclidean space are always at least $60^\circ$.
The Minkowski analogue, observed by Cieslik \cite{Cieslik}, is as follows.
\begin{lemma}\label{lemma1}
If $ba$ and $bc$ are two edges in an MST in a Minkowski space, then $\size{\myangle abc}\geq 1$.
\end{lemma}
\begin{proof}
Without loss of generality we may assume that $\norm{a-b}\geq\norm{c-b}$ (otherwise interchange $a\leftrightarrow c$).
Let
\[ d=b+\frac{\norm{c-b}}{\norm{a-b}}(a-b).\]
Then $d$ is on the segment $ab$, and $\norm{b-d}=\norm{b-c}$.
Since $ba$ and $bc$ are edges of the MST, $\norm{a-c}\geq\norm{a-b}$.
On the other hand, by the triangle inequality, $\norm{a-c}\leq\norm{a-d}+\norm{d-c}$, and we obtain $\norm{d-c}\geq\norm{a-b}-\norm{a-d}=\norm{b-d}$.
Substituting the definition of $d$ into this inequality and dividing by $\norm{c-b}$ we obtain $\size{\myangle abc}\geq 1$.
\end{proof}

Choose any $\epsi < \min\{\norm{x-y}:x,y\in S,x\neq y\}$.
An \emph{$\epsi$-perturbation} of $S$ is any set of points $S'$ in $X^d$ for which there exists a bijection $f:S\to S'; x\mapsto x'$ such that $\abs{\norm{x-y}-\norm{x'-y'}}\leq\epsi$ for all $x,y\in S$.
By the choice of $\epsi$ we may identify the $\epsi$-perturbation with the function $f$.
The set $\CP=\prod_{x\in S} B(x,\epsi)$ of all $\epsi$-perturbations of $S$ can be given the metric $d(f_1,f_2)=\max\{\norm{f_1(x)-f_2(x)}:x\in S\}$.
Then $\CP$ is a complete metric space.
We now show that any $S$ has an $\epsi$-perturbation $S'$ such that no angle $\myangle a'b'c'$ determined by $a',b',c'\in S'$ has size $1$.
Such an $\epsi$-perturbation will be any point in $\CP$ outside the closed set $\CS_{\myangle abc}=\{S'\in\CP: \size{\myangle a'b'c'}=1\}$, for any distinct $a,b,c\in S$.
By the Baire category theorem \cite[\S12]{Simmons} it is sufficient to prove that each $\CS_{\myangle abc}$ is nowhere dense in $\CP$.
This follows from the following lemma in planar Minkowski geometry.

\begin{lemma}
Let $\myangle abc$ be any angle in a Minkowski plane $X^2$ such that $\size{\myangle abc}=1$.
Then for any $\epsi>0$ there exists an angle $\myangle a'b'c'$ with $\norm{a-a'}\leq\epsi$, $\norm{b-b'}\leq\epsi$, $\norm{c-c'}\leq\epsi$ and $\size{\myangle a'b'c'}\neq 1$.
\end{lemma}
\begin{proof}
Without loss of generality we may assume that $b=o$.
Let $x=\frac{1}{\norm{a}}a$ and $y=\frac{1}{\norm{c}}c$.
If for all $a'\in B(a,\epsi)$ and $c'\in B(c,\epsi)$ we still have $\size{\myangle a'bc'}=1$, then all chords of the unit ball parallel to $xy$ and sufficiently close to the chord $xy$ have length $1$.
Since the unit ball is convex, this is only possible if $x$ and $y$ are both contained in two parallel segments on the boundary of the unit ball.
However, such parallel segments are either on the same line or on two different lines at distance $2$ from each other.
Both cases give a contradiction.
\end{proof}
The proof of the above lemma also gives that $\{S'\in\CP: \size{\myangle a'b'c'}=\lambda\}$ is nowhere dense in $\CP$ for any $0<\lambda<2$.
For $\lambda=2$ this set is not necessarily nowhere dense if the norm is not strictly convex.

Up to now we have shown that for any sufficiently small $\epsi>0$ there exists an $\epsi$-perturbation $S_\epsi$ of $S$ such that no angle in $S_\epsi$ has size $1$.
Consider now an MST of $S_\epsi$.
Let $o$ be a point in $S_\epsi$ which has largest degree, and let its neighbours be $p_1,\dots,p_k$.
By Lemma~\ref{lemma1}, $\size{\myangle p_i o p_j}\geq 1$, and by choice of $S_\epsi$, $\size{\myangle p_i o p_j}\neq 1$ for any $i\neq j$.
It follows that if we let $x_i=\norm{p_i}^{-1}p_i$, then $x_1,\dots,x_k$ are unit vectors with $\norm{x_i-x_j}>1$ for all distinct $i,j$.
Therefore, $k\leq H_s(B)$.
If we now let $\epsi=1/n$ for sufficiently large $n$, we obtain a sequence of MSTs, each with maximum degree bounded above by $H_s(B)$.
Since there are only finitely many trees on a finite set of points, there is a subsequence with the same tree structure.
This subsequence converges to a tree on $S$ which, by continuity of the norm, is an MST of $S$.
We have found an MST $T$ of $S$ with $\Delta(T)\leq H_s(B)$.
This shows that $\Delta^-(X^d)\leq H_s(B)$, and finishes the proof of the theorem.


\end{document}